\documentclass[english]{article}
\usepackage[T1]{fontenc}
\usepackage[latin9]{inputenc}
\usepackage{float}
\usepackage{amsmath}
\usepackage{amsthm}
\usepackage{amssymb}
\usepackage{graphicx}

\makeatletter

\providecommand{\tabularnewline}{\\}

\numberwithin{equation}{section}
\numberwithin{figure}{section}
\numberwithin{table}{section}
\newcommand{\lyxaddress}[1]{
	\par {\raggedright #1
	\vspace{1.4em}
	\noindent\par}
}

\makeatother

\usepackage{babel}
\begin{document}
\title{Ultradiscrete two-variable Oregonator}
\author{Yuhei KASHIWATATE}
\maketitle

\lyxaddress{\begin{center}
Graduate School of Science and Engineering, Yamagata University,
Kojirakawa-machi 1-4-12, Yamagata 990-8560, Japan\\
E-mail: s181104m@st.yamagata-u.ac.jp
\par\end{center}}
\begin{abstract}
Ultradiscretization is a limiting procedure transforming a given differential/difference
equation into a ultradiscrete equation. Ultradiscrete equations are
expressed by addition, subtraction and/or max. The procedure is expected
to preserve the essential properties of the original equations. As
a method of ultradiscretization, there is ``tropical discretization''
proposed by M. Murata. In this paper, we shall modify it, and derive
a ultradiscrete equation from the continuous model of the BZ reaction.
The derived equation generates a cellular automaton by restricting
the values of the parameters, which is equivalent to one of those
introduced by D. Takahashi, A. Shida, and M. Usami. By setting appropriate
initial values, we can obtain the typical patterns of the BZ reaction.
Furthermore, we consider the equation without diffusion effect and
derive the explicit solutions. As a result, the solutions corresponding
to the limit cycle (oscillation) appearing in the continuous model
will be found.
\end{abstract}

\section{Introduction}

Ultradiscretization is a limiting procedure transforming a given differential/diffe-
rence equation into a ultradiscrete equation. Ultradiscrete equations
are based on the max-plus algebra, which defines the summation by
max and the product by +, and often expressed in simple forms. It
is expected that this procedure preserves the essential properties
of the original equations. In the reaction-diffusion system, both
continuous models using partial differential equations and mathematical
models using cellular automatons have been studied. However, the direct
correspondence between them is not clear. Indeed, it was difficult
to obtain max-plus equations expressing the BZ reaction directly from
the system of partial differential equations. For example, in their
paper \cite{key-3}, D. Takahashi, A. Shida, and M. Usami created
max-plus equations for the BZ reaction simply by comparing its typical
nature.

As a method of ultradiscretization, there is ``tropical discretization''
proposed by M. Murata \cite{key-1,key-2}. In this paper, we shall
modify it, and derive a ultradiscrete equation from the continuous
model of the BZ reaction. Ultradiscretization of two-dimensional space
diffusion terms will be the mean of values of neighbouring four points
and center. In Murata's method, the value of center is not used. Another
difference is which terms are shifted in time. The derived equation
generates a cellular automaton by restricting the values of the parameters,
which is equivalent to one of those introduced in the paper \cite{key-3}.
By setting appropriate initial values, we obtain the typical patterns
of the BZ reaction, ring, target and spiral, which are mentioned in
the book \cite{key-4}.

Furthermore, we consider the equation without diffusion effect and
derive the explicit solutions. As a result, the solutions corresponding
to the limit cycle (oscillation) appearing in the continuous model
will be found.

\section{The two-variable Oregonator}

The two-variable Oregonator in two-dimensional space is given by
\begin{equation}
\left\{ \begin{aligned}\frac{\partial u}{\partial t} & =D_{u}\left(\frac{\partial^{2}u}{\partial x^{2}}+\frac{\partial^{2}u}{\partial y^{2}}\right)+a\left\{ u(1-u)-\frac{fv(u-q)}{u+q}\right\} ,\\
\frac{\partial v}{\partial t} & =D_{v}\left(\frac{\partial^{2}v}{\partial x^{2}}+\frac{\partial^{2}v}{\partial y^{2}}\right)+u-v,
\end{aligned}
\right.\label{eq:-4}
\end{equation}
where $u=u(t,x,y)$ is the activity factor and $v=v(t,x,y)$ is the
inhibitory factor. The constants $D_{u}$ and $D_{v}$ are the diffusion
coefficients of $u$ and $v$, respectively. We consider $(x,y)\in\mathbb{R}^{2},\:t\geq0$
and the constants $a,f,q$ satisfy $a\sim0.25\times10^{2},\:1<f<2,\:q\sim8\times10^{-4}$.
This system is known as a model to explain the pattern dynamics of
the BZ reaction. The solutions of this system represent spatial patterns.
Changing initial values and the values of paremeters, we can observe
various patterns. 

\subsection{Discretization}

In this section, we discretize the eq.(\ref{eq:-4}) by a method similar
to Murata's \cite{key-1}. First we consider a discretization of the
following system of partial differential equations:
\begin{equation}
\begin{cases}
{\displaystyle {\displaystyle \frac{\partial u}{\partial t}}=D_{u}\left({\displaystyle \frac{\partial^{2}u}{\partial x^{2}}+\frac{\partial^{2}u}{\partial y^{2}}}\right)} & (D_{u}>0),\\
{\displaystyle {\displaystyle \frac{\partial v}{\partial t}}=D_{v}\left({\displaystyle \frac{\partial^{2}v}{\partial x^{2}}+\frac{\partial^{2}v}{\partial y^{2}}}\right)} & (D_{v}>0).
\end{cases}\label{eq:}
\end{equation}
The following system of difference equations is that.
\begin{equation}
\left\{ \begin{aligned}u_{n+1}^{j,k} & =\frac{1}{5}(u_{n}^{j,k}+u_{n}^{j-\alpha,k}+u_{n}^{j+\alpha,k}+u_{n}^{j,k-\alpha}+u_{n}^{j,k+\alpha})=m_{\alpha}(u_{n}^{j,k}),\\
v_{n+1}^{j,k} & =\frac{1}{5}(v_{n}^{j,k}+v_{n}^{j-\beta,k}+v_{n}^{j+\beta,k}+v_{n}^{j,k-\beta}+v_{n}^{j,k+\beta})=m_{\beta}(v_{n}^{j,k}).
\end{aligned}
\right.\label{eq:-1}
\end{equation}
Indeed, if we put $u_{n}^{j,k}=u(t,x,y)$ with $t=n\Delta t,\:x=j\Delta x,\:y=k\Delta y\:(\Delta x=\Delta y)$,
we find
\begin{align*}
u_{n+1}^{j,k} & =u(t+\Delta t,x,y),\\
u_{n}^{j\pm\alpha,k} & =u(t,x\pm\alpha\Delta x,y),\\
u_{n}^{j,k\pm\alpha} & =u(t,x,y\pm\alpha\Delta y).
\end{align*}
By the taylor expansions at $(t,x,y)$, we get
\begin{align*}
u(t,x,y) & =u,\\
u(t+\Delta t,x,y) & =u+u_{t}\Delta t+\frac{1}{2}u_{tt}\Delta t^{2}+\cdots,\\
u(t,x\pm\alpha\Delta x,y) & =u\pm u_{x}\alpha\Delta x+\frac{1}{2}u_{xx}\alpha^{2}\Delta x^{2}\pm\cdots,\\
u(t,x,y\pm\alpha\Delta y) & =u\pm u_{y}\alpha\Delta y+\frac{1}{2}u_{yy}\alpha^{2}\Delta y^{2}\pm\cdots.
\end{align*}
Substituting them into the eq.(\ref{eq:-1}), we get
\begin{align*}
u+u_{t}\Delta t+\frac{1}{2}u_{tt}\Delta t^{2}+\cdots=\frac{1}{5} & \left\{ u+\left(u+u_{x}\alpha\Delta x+\frac{1}{2}u_{xx}\alpha^{2}\Delta x^{2}+\cdots\right)\right.\\
 & \left.+\left(u-u_{x}\alpha\Delta x+\frac{1}{2}u_{xx}\alpha^{2}\Delta x^{2}-\cdots\right)\right.\\
 & \left.+\left(u+u_{y}\alpha\Delta y+\frac{1}{2}u_{yy}\alpha^{2}\Delta y^{2}+\cdots\right)\right.\\
 & \left.+\left(u-u_{y}\alpha\Delta y+\frac{1}{2}u_{yy}\alpha^{2}\Delta y^{2}-\cdots\right)\right\} ,
\end{align*}
and thus
\begin{align}
u_{t}\Delta t+\cdots= & \frac{1}{5}u_{xx}\alpha^{2}\Delta x^{2}+\frac{1}{5}u_{yy}\alpha^{2}\Delta y^{2}+\cdots,\nonumber \\
u_{t}+\cdots= & \frac{\alpha^{2}\Delta x^{2}}{5\Delta t}u_{xx}+\frac{\alpha^{2}\Delta y^{2}}{5\Delta t}u_{yy}+\cdots,\label{eq:-10}
\end{align}
where $D=\alpha^{2}\Delta x^{2}/5\Delta t=\alpha^{2}\Delta y^{2}/5\Delta t$.
Taking the limit $\Delta t,\Delta x,\Delta y\rightarrow+0$ without
change of $D$, we obtain the first equation of the system eq.(\ref{eq:}).
Thus, the eq.(\ref{eq:-1}) can be regarded as a discretization of
the eq.(\ref{eq:}).

Furthermore, we consider a discretization of the following syestem
of ordinary differential equations:
\begin{equation}
\left\{ \begin{aligned}\frac{du}{dt} & =a\left\{ u(1-u)-\frac{fv(u-q)}{u+q}\right\} ,\\
\frac{dv}{dt} & =u-v.
\end{aligned}
\right.\label{eq:-2}
\end{equation}
We shall show that the following syestem of difference equations is
the required
\begin{equation}
\left\{ \begin{aligned}u_{n+1} & =\frac{\varepsilon^{-1}u_{n}+au_{n}+{\displaystyle \frac{afqv_{n}}{u_{n}+q}}}{\varepsilon^{-1}+au_{n}+{\displaystyle \frac{afv_{n}}{u_{n}+q}}},\\
v_{n+1} & =\frac{\varepsilon^{-1}v_{n}+u_{n}}{\varepsilon^{-1}+1},
\end{aligned}
\right.\label{eq:-3}
\end{equation}
where $n\in\mathbb{Z}_{\geq0},\:\varepsilon>0$. The method we adopt
here is the same as that in the papers \cite{key-1,key-2}.

Putting $u_{n}=u(t),\:v_{n}=v(t),\:t=\varepsilon n$, we find
\begin{equation}
\left\{ \begin{aligned}\frac{u(t+\varepsilon)-u(t)}{\varepsilon} & =a\left\{ u(t)(1-u(t))-\frac{fv(t)(u(t)-q)}{u(t)+q}\right\} +O(\varepsilon),\\
\frac{v(t+\varepsilon)-v(t)}{\varepsilon} & =u(t)-v(t)+O(\varepsilon).
\end{aligned}
\right.\label{eq:-5}
\end{equation}
Taking the limit $\varepsilon\rightarrow+0$, we obtain the system
of differential equations(\ref{eq:-2}). Thus, the eq.(\ref{eq:-3})
can be regarded as a discretization of the eq.(\ref{eq:-2}). Using
the eq.(\ref{eq:-1}) and the eq.(\ref{eq:-3}), we will find the
system of difference equations,
\begin{equation}
\left\{ \begin{aligned}u_{n+1}^{j,k} & =\frac{\epsilon^{-1}m_{\alpha}(u_{n}^{j,k})+am_{\alpha}(u_{n}^{j,k})+{\displaystyle \frac{afqv_{n}^{j,k}}{m_{\alpha}(u_{n}^{j,k})+q}}}{\epsilon^{-1}+am_{\alpha}(u_{n}^{j,k})+{\displaystyle \frac{afv_{n}^{j,k}}{m_{\alpha}(u_{n}^{j,k})+q}}},\\
v_{n+1}^{j,k} & =\frac{\epsilon^{-1}m_{\beta}(v_{n}^{j,k})+m_{\beta}(u_{n}^{j,k})}{\epsilon^{-1}+1}.
\end{aligned}
\right.\label{eq:-6}
\end{equation}
This can be rewritten as
\begin{equation}
\left\{ \begin{aligned}\frac{u_{n+1}^{j,k}-u_{n}^{j,k}}{\varepsilon} & =\frac{m_{\alpha}(u_{n}^{j,k})-u_{n}^{j,k}}{\varepsilon}+a\left\{ m_{\alpha}(u_{n}^{j,k})(1-u_{n}^{j,k})-\frac{fv_{n}^{j,k}(u_{n}^{j,k}-q)}{m_{\alpha}(u_{n}^{j,k})+q}\right\} +O(\varepsilon),\\
\frac{v_{n+1}^{j,k}-v_{n}^{j,k}}{\varepsilon} & =\frac{m_{\beta}(v_{n}^{j,k})-v_{n}^{j,k}}{\varepsilon}+m_{\beta}(u_{n}^{j,k})-v_{n}^{j,k}+O(\varepsilon),
\end{aligned}
\right.\label{eq:-11}
\end{equation}
and thus can be regarded as a discretization of the eq.(\ref{eq:-4}).

\subsection{Ultradiscretization}

In this section, we shall ultradiscretize the eq.(\ref{eq:-6}) and
investigate the solutions. Let
\begin{equation}
\begin{array}{ccc}
u_{n}^{j,k}=\exp(U_{n}^{j,k}/\lambda), & v_{n}^{j,k}=\exp(V_{n}^{j,k}/\lambda), & \epsilon=\exp(E/\lambda),\\
a=\exp(A/\lambda), & f=\exp(F/\lambda), & q=\exp(Q/\lambda),
\end{array}\label{eq:2,10}
\end{equation}
and take the limit $\lambda\rightarrow+0$. Operations such as
\[
\begin{aligned}\lim_{\lambda\rightarrow+0}\lambda\log\left(e^{A/\lambda}+e^{B/\lambda}\right) & =\max(A,B),\\
\lim_{\lambda\rightarrow+0}\lambda\log\left(e^{A/\lambda}\cdotp e^{B/\lambda}\right) & =A+B,\\
\lim_{\lambda\rightarrow+0}\lambda\log\left(e^{A/\lambda}\text{/}e^{B/\lambda}\right) & =A-B,
\end{aligned}
\]
perform here. Therefore, the eq.(\ref{eq:-6}) is transformed into
\begin{equation}
\left\{ \begin{aligned}U_{n+1}^{j,k}= & \max\{M_{\alpha}(U_{n}^{j,k})-E,A+M_{\alpha}(U_{n}^{j,k}),A+F+Q+V_{n}^{j,k}-\max(M_{\alpha}(U_{n}^{j,k}),Q)\}\\
 & -\max\{-E,A+M_{\alpha}(U_{n}^{j,k}),A+F+V_{n}^{j,k}-\max(M_{\alpha}(U_{n}^{j,k}),Q)\},\\
V_{n+1}^{j,k}= & \max\{M_{\beta}(V_{n}^{j,k})-E,M_{\beta}(U_{n}^{j,k})\}-\max\{-E,0\},
\end{aligned}
\right.\label{eq:-7}
\end{equation}
where
\[
\left\{ \begin{aligned}M_{\alpha}(U_{n}^{j,k})= & \max(U_{n}^{j,k},U_{n}^{j-\alpha,k},U_{n}^{j+\alpha,k},U_{n}^{j,k-\alpha},U_{n}^{j,k+\alpha}),\\
M_{\beta}(V_{n}^{j,k})= & \max(V_{n}^{j,k},V_{n}^{j-\beta,k},V_{n}^{j+\beta,k},V_{n}^{j,k-\beta},V_{n}^{j,k+\beta}),
\end{aligned}
\right.
\]
which is an ultradiscretization of the eq.(\ref{eq:-1}).

Taking the limit $E\rightarrow+\infty$, we get
\begin{equation}
\left\{ \begin{aligned}U_{n+1}^{j,k}= & \max\{M_{\alpha}(U_{n}^{j,k}),F+Q+V_{n}^{j,k}-\max(M_{\alpha}(U_{n}^{j,k}),Q)\}\\
 & -\max\{M_{\alpha}(U_{n}^{j,k}),F+V_{n}^{j,k}-\max(M_{\alpha}(U_{n}^{j,k}),Q)\},\\
V_{n+1}^{j,k}= & M_{\beta}(U_{n}^{j,k}).
\end{aligned}
\right.\label{eq:-8}
\end{equation}
and the following single ultradiscrete equation:
\begin{align}
U_{n+1}^{j,k}= & \max\{M_{\alpha}(U_{n}^{j,k}),F+Q+M_{\beta}(U_{n-1}^{j,k})-\max(M_{\alpha}(U_{n}^{j,k}),Q)\}\nonumber \\
 & -\max\{M_{\alpha}(U_{n}^{j,k}),F+M_{\beta}(U_{n-1}^{j,k})-\max(M_{\alpha}(U_{n}^{j,k}),Q)\}.\label{eq:2.13}
\end{align}

\section{Cellular automaton}

Considering the values of $f$ and $q$ in the eq.(\ref{eq:-4}),
we suppose $Q<0<F$ in the eq.(\ref{eq:2.13}). We take $Q=-1,\:F=1$,
and restrict initial values to $U_{0}^{j,k},U_{1}^{j,k}\in\{-1,0\}$
to generate a cellulor automaton. Under this condition, the values
of $U_{n}^{j,k}$ are restricted to $\{-1,0\}$. To shift the values
to $\{0,1\}$, Let
\begin{equation}
W_{n}^{j,k}=U_{n}^{j,k}-Q,\label{eq:2.14}
\end{equation}
and substitute it into the eq.(\ref{eq:2.13}). We get
\begin{align}
W_{n+1}^{j,k}= & \max\{M_{\alpha}(W_{n}^{j,k}),F+M_{\beta}(W_{n-1}^{j,k})-M_{\alpha}(W_{n}^{j,k})\}\nonumber \\
 & -\max\{M_{\alpha}(W_{n}^{j,k})+Q,F+M_{\beta}(W_{n-1}^{j,k})-M_{\alpha}(W_{n}^{j,k})\}.\label{eq:2.17-1}
\end{align}
Since we find
\[
\begin{aligned}M_{\alpha}(W_{n}^{j,k}) & \in\{0,1\},\\
F+M_{\beta}(W_{n-1}^{j,k})-M_{\alpha}(W_{n}^{j,k}) & \in\{0,1,2\},\\
M_{\alpha}(W_{n}^{j,k})+Q & \in\{-1,0\},
\end{aligned}
\]
it follows that
\[
F+M_{\beta}(W_{n-1}^{j,k})-M_{\alpha}(W_{n}^{j,k})\geq M_{\alpha}(W_{n}^{j,k})+Q.
\]
Therefore, we obtain the following simple equation:
\begin{equation}
W_{n+1}^{j,k}=\max\{2M_{\alpha}(W_{n}^{j,k})-M_{\beta}(W_{n-1}^{j,k})-F,0\}.\label{eq:2.18}
\end{equation}
The rule of time evolution is in the Table 3.1:
\begin{table}[H]
\centering{}%
\begin{tabular}{|c||c|c|c|c|}
\hline 
$M_{\alpha}(W_{n}^{j,k}),M_{\beta}(W_{n-1}^{j,k})$ & $0,0$ & $0,1$ & $1,0$ & $1,1$\tabularnewline
\hline 
$W_{n+1}^{j,k}$ & $0$ & $0$ & $1$ & $0$\tabularnewline
\hline 
\end{tabular}\caption{$W_{n+1}^{j,k}$\label{Table 3.1}}
\end{table}

We shall compare it with one of Takahashi, Shida and Usami's max-plus
equations \cite{key-3}:
\begin{equation}
Y_{n+1}^{j,k}=\max(Y_{n}^{j,k},Y_{n}^{j-1,k},Y_{n}^{j+1,k},Y_{n}^{j,k-1},Y_{n}^{j,k+1},Y_{n-1}^{j,k})-Y_{n-1}^{j,k}.\label{eq:2.15}
\end{equation}
We recall thet this equation was not associated with the differential
equation directly. Using our terminology, we obtain the following
form:
\begin{equation}
\begin{array}{cc}
Y_{n+1}^{j,k} & =\max(M_{1}(Y_{n}^{j,k})-Y_{n-1}^{j,k},0)\end{array}.\label{eq:2.16}
\end{equation}
The rule of time evolution is in the Table 3.2:
\begin{table}[H]
\centering{}%
\begin{tabular}{|c||c|c|c|c|}
\hline 
$M_{1}(Y_{n}^{j,k}),Y_{n-1}^{j,k}$ & $0,0$ & $0,1$ & $1,0$ & $1,1$\tabularnewline
\hline 
$Y_{n+1}^{j,k}$ & $0$ & $0$ & $1$ & $0$\tabularnewline
\hline 
\end{tabular}\caption{$Y_{n+1}^{j,k}$\label{Table 3.2}}
\end{table}
In the case of $(\alpha,\beta)=(1,0)$ the rules are the same. Therefore,
we have found a connection between the cellular automaton and the
two-variable Oregonator by tropical discretization. For reader's convenience,
we shall introduce examples of simulation, because there is no example
for this cellular automaton in the paper \cite{key-3}. The following
are introduced in the text-book \cite{key-4}.

The Figure 3.1 shows a `single ring' pattern. From the center, a
square-shaped wave with value 1 spreads outwards.
\begin{figure}[H]
\centering{}\includegraphics[scale=0.4]{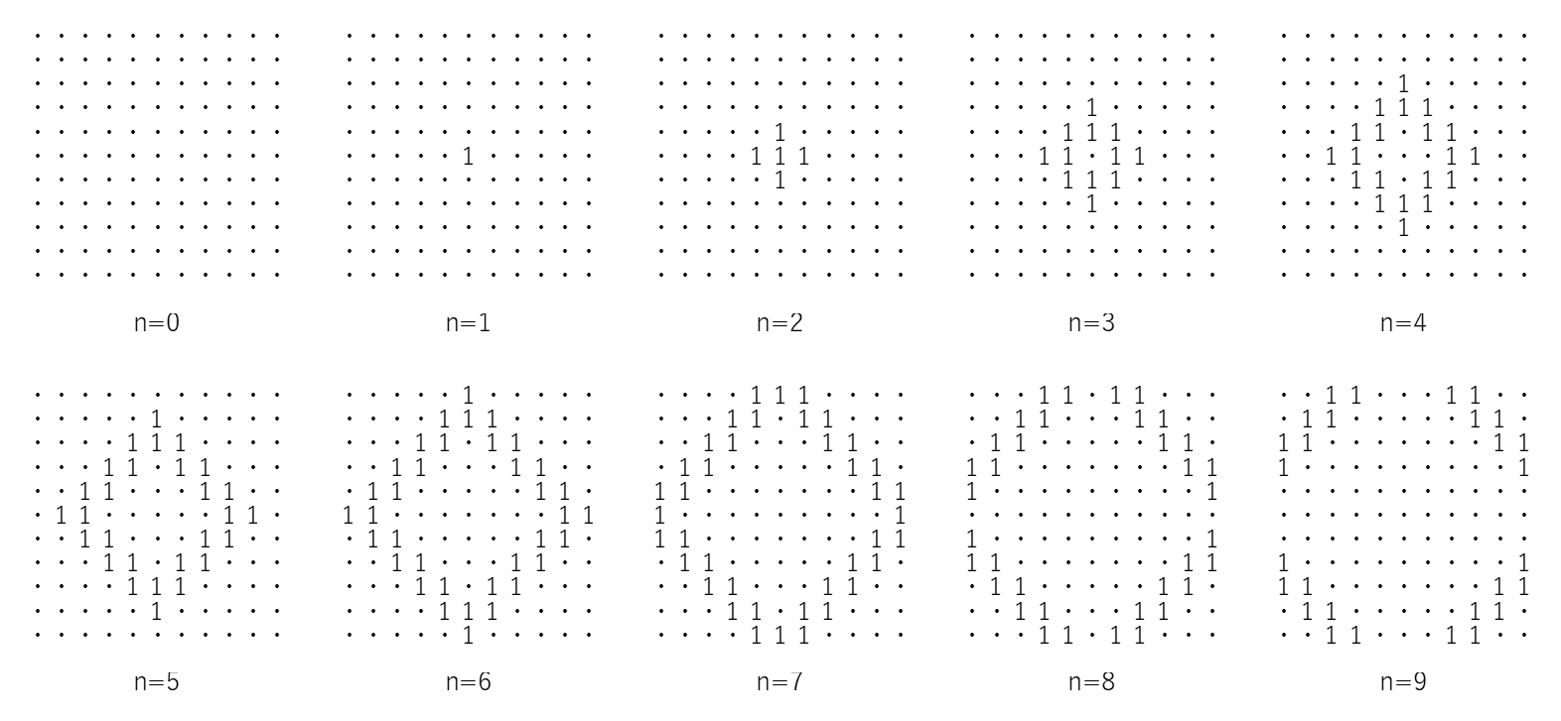}\caption{Single ring pattern.\label{Figure 3.1}}
\end{figure}

The Figure 3.2 shows a process to form a stable `target' pattern.
At the center, the value changes periodically as 1,1,0,0 and square-shaped
waves appear and spread outwards repeatedly with period 4.
\begin{figure}[H]
\centering{}\includegraphics[scale=0.4]{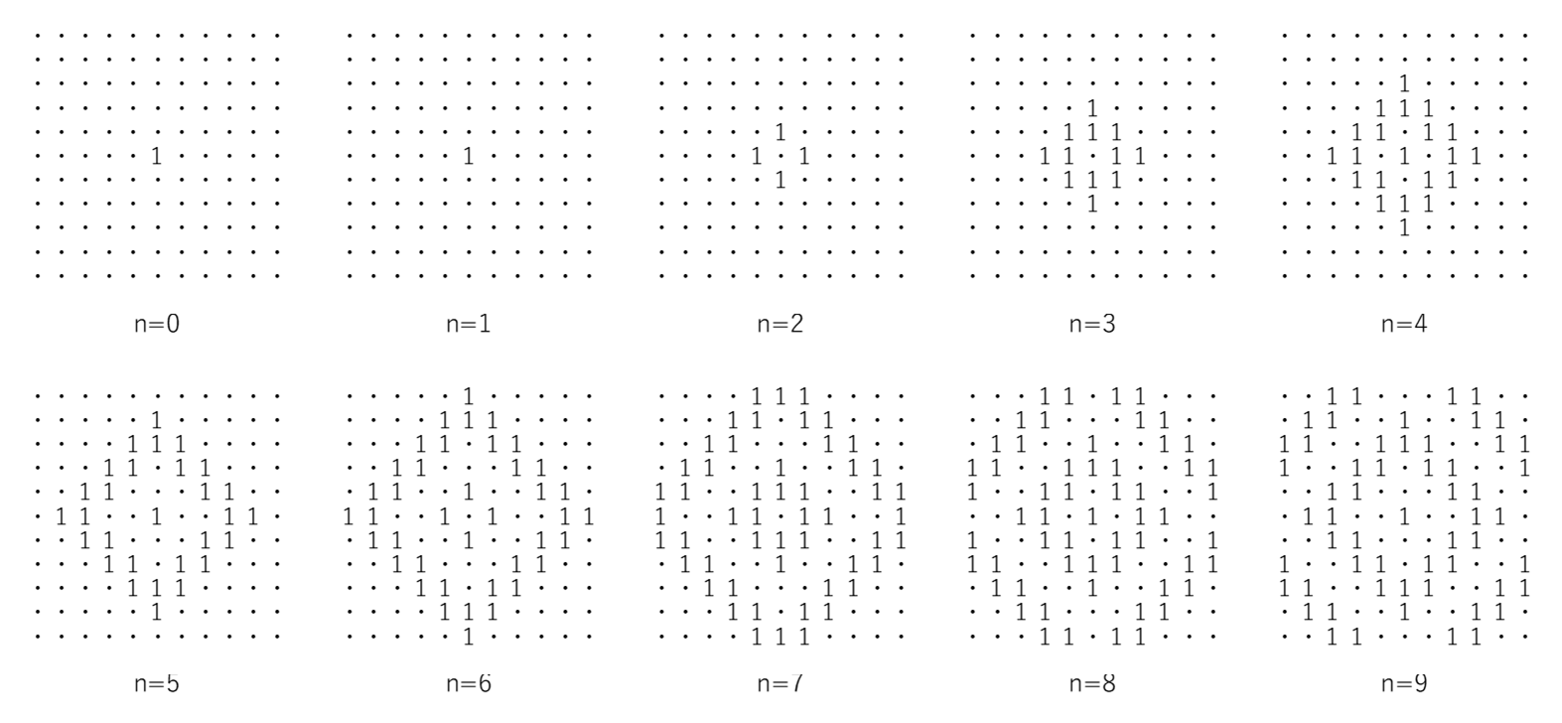}\caption{Target pattern.\label{Figure 3.2}}
\end{figure}

The Figure 3.3 shows a process to form a stable `spiral' pattern.
By an horizontal line of value 1, a spirals appear from its end points.
Thereafter, the spirals spread through the whole space and each rotates
by 90 degrees per unit time. Also, the spirals collide and disappear
without any other interaction.
\begin{figure}[H]
\begin{centering}
\includegraphics[scale=0.3]{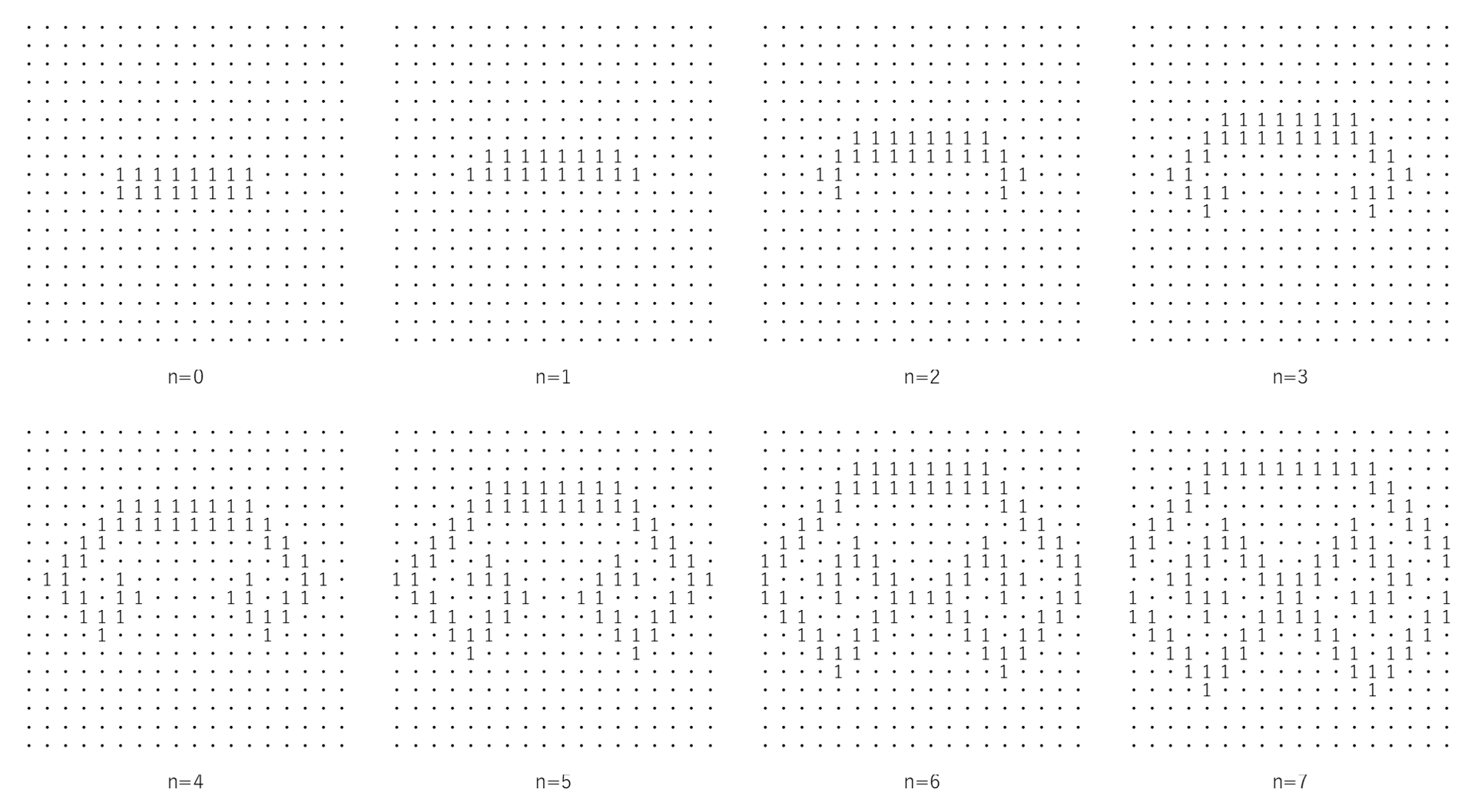}
\par\end{centering}
\caption{Spiral pattern.\label{Figure 3.3}}
\end{figure}
The behavior of these solutions is similar to that of the BZ reaction
patterns.

\section{Ultradiscrete two-variable Oregonator}

\subsection{Equilibrium points}

The eq.(\ref{eq:2.13}) with $(\alpha,\beta)=(0,0)$ is the following
second order ordinary difference equation:
\begin{align}
U_{n+1}= & \max\{U_{n},F+Q+U_{n-1}-\max(U_{n},Q)\}\nonumber \\
 & -\max\{U_{n},F+U_{n-1}-\max(U_{n},Q)\},\label{eq:-12}
\end{align}
which can be regarded as the equation without diffusion effect. We
consider the equilibrium points of the eq.(\ref{eq:-12}).

(I) If $U_{n}\geq Q$, we get
\begin{equation}
\bar{U}=\begin{cases}
0 & (F+Q\leq0\;\&\;F\leq0),\\
Q & (F\geq0\;\&\;F\geq Q).
\end{cases}\label{eq:-21}
\end{equation}

(II) If $U_{n}\leq Q$, we get
\begin{equation}
\bar{U}=\begin{cases}
0 & (F\leq Q\;\&\;F\leq0),\\
F & (0<F<Q),\\
Q & (F\geq0\;\&\;F\geq Q),
\end{cases}\label{eq:-22}
\end{equation}
where, $0,\:Q$ are stable equilibrium points, and $F$ is an unstable
equilibrium point.

\subsection{Periodic solution}

In this section, we suppose $0<F<Q$ and $F,\:Q,\:U_{0},\:U_{1}\in\mathbb{Z}$
in the eq.(\ref{eq:-12}), which implies $U_{n}\in\mathbb{Z}$. In
this case, the stable equilibrium points disappear. This situaton
is very similar to that of the original differential equations in
which the limit cycle (oscillation) appears.

We shall obtain the solution of the eq.(\ref{eq:-12}).

(I) When $U_{n}\geq Q$, we get
\begin{equation}
U_{n+1}=\begin{cases}
0 & (F+Q+U_{n-1}-2U_{n}\leq0),\\
F+Q+U_{n-1}-2U_{n} & (F+U_{n-1}-2U_{n}<0<F+Q+U_{n-1}-2U_{n}),\\
Q & (F+U_{n-1}-2U_{n}\geq0),
\end{cases}\label{eq:-23}
\end{equation}
which implies $U_{n+1}\leq Q$. Thus, the next step is in the following
case (II).

(II) When $U_{n}\leq Q$, we get 
\begin{equation}
U_{n+1}=\begin{cases}
0 & (F+U_{n-1}-U_{n}\leq0),\\
F+U_{n-1}-U_{n} & (F-Q+U_{n-1}-U_{n}<0<F+U_{n-1}-U_{n}),\\
Q & (F-Q+U_{n-1}-U_{n}\geq0).
\end{cases}\label{eq:-24}
\end{equation}
It follows that $U_{n+1}\leq Q$. Thus, once $U_{n}$ satisfies the
case (II), it keeps in the case (II) all the time.

Therefore, we only consider solutions satisfying the case (II). For
brevity, let $\Psi_{n}=U_{n-1}-U_{n}$. Eq.(\ref{eq:-24}) can be
rewritten as
\begin{equation}
U_{n+1}=\begin{cases}
0 & (\Psi_{n}\leq-F),\\
F+\Psi_{n} & (-F<\Psi_{n}<Q-F),\\
Q & (\Psi_{n}\geq Q-F).
\end{cases}\label{eq:-25}
\end{equation}

Next, we put $I_{1}=(-\infty,-F],\:I_{2}=(-F,Q-F),\:I_{3}=[Q-F,\infty)$.
The combination of $(\Psi_{n},\Psi_{n+1})$ can be divided into $3\times3$
cases.

(1-1) In the case $(\Psi_{n},\Psi_{n+1})\in I_{1}\times I_{1}$, we
get
\begin{align*}
(U_{n+1},\Psi_{n+1}) & =(0,U_{n}),\\
(U_{n+2},\Psi_{n+2}) & =(0,0),
\end{align*}
and thus $(\Psi_{n+1},\Psi_{n+2})\in I_{1}\times I_{2}$. Therefore,
the state changes from (1-1) to (1-2).

(1-2) In the case $(\Psi_{n},\Psi_{n+1})\in I_{1}\times I_{2}$, it
follows that $(\Psi_{n+1},\Psi_{n+2})\in I_{2}\times I_{i}$ for some
$i$.

(1-3) In the case $(\Psi_{n},\Psi_{n+1})\in I_{1}\times I_{3}$, we
get
\begin{align*}
(U_{n+1},\Psi_{n+1}) & =(0,U_{n}),\\
(U_{n+2},\Psi_{n+2}) & =(Q,-Q),\\
(U_{n+3},\Psi_{n+3}) & =(0,Q),\\
(U_{n+4},\Psi_{n+4}) & =(Q,-Q),\\
 & \vdots
\end{align*}
Therefore, $U_{n}$ alternates between $0$ and $Q$ from a certain
time.

(2-1) In the case $(\Psi_{n},\Psi_{n+1})\in I_{2}\times I_{1}$, we
get
\begin{align*}
(U_{n+1},\Psi_{n+1}) & =(F+\Psi_{n},U_{n}-(F+\Psi_{n})),\\
(U_{n+2},\Psi_{n+2}) & =(0,F+\Psi_{n}),
\end{align*}
where
\[
\begin{cases}
\Psi_{n+2}\in I_{3} & (\Psi_{n}\geq Q-2F),\\
\Psi_{n+2}\in I_{2} & (\Psi_{n}<Q-2F).
\end{cases}
\]

\,\,(2-1-a) In the case $\Psi_{n}\geq Q-2F$, we have $(\Psi_{n+1},\Psi_{n+2})\in I_{1}\times I_{3}$.

\,\,(2-1-b) In the case $\Psi_{n}<Q-2F$, we get
\begin{align*}
(U_{n+3},\Psi_{n+3}) & =(2F+\Psi_{n},-2F-\Psi_{n}),\\
(U_{n+4},\Psi_{n+4}) & =(0,2F+\Psi_{n}),
\end{align*}
where
\[
\begin{cases}
\Psi_{n+4}\in I_{3} & (\Psi_{n}\geq Q-3F),\\
\Psi_{n+4}\in I_{2} & (\Psi_{n}<Q-3F).
\end{cases}
\]

\,\,(2-1-c) In the case $\Psi_{n}\geq Q-3F$, we have $(\Psi_{n+3},\Psi_{n+4})\in I_{1}\times I_{3}$.

\,\,(2-1-d) In the case $\Psi_{n}<Q-3F$, repeat the above operation.

If the above discussion did not stop, $\Psi_{n+i}\notin I_{3}$ would
be satisfied all the time. That means $\Psi_{n}<Q-mF\:\:(m\geq2)$,
which contradicts $-F<\Psi_{n}<Q-F$. Therefore, $(\Psi_{n'},\Psi_{n'+1})\in I_{1}\times I_{3}$
is satisfied for a certain $n'$.

(2-2) In the case $(\Psi_{n},\Psi_{n+1})\in I_{2}\times I_{2}$, it
follows that $(\Psi_{n+1},\Psi_{n+2})\in I_{2}\times I_{i}$ for some
$i$. If we consider the case of continuing to meet $(\Psi_{n},\Psi_{n+1})\in I_{2}\times I_{2}$,
we get
\begin{equation}
\Psi_{n+2}=\Psi_{n}-\Psi_{n+1}.\label{eq:-26}
\end{equation}
The general solution is
\begin{align}
\Psi_{n}= & (-1)^{n-1}c_{1}\left(\frac{1+\sqrt{5}}{2}\right)^{n-1}+c_{2}\left(\frac{-1+\sqrt{5}}{2}\right)^{n-1},\label{eq:-9}
\end{align}
where 
\[
c_{1}=\Psi_{2}-\left(\frac{-1+\sqrt{5}}{2}\right)\Psi_{1},\:c_{2}=\Psi_{2}+\left(\frac{1+\sqrt{5}}{2}\right)\Psi_{1}.
\]
We easily find $c_{1}=0\Leftrightarrow\Psi_{1}=\Psi_{2}=0$ by $U_{n}\in\mathbb{Z}$.

\,\,(2-2-a) In the case $(\Psi_{1},\Psi_{2})\neq(0,0)$, taking
the limit $n\rightarrow\infty$, absolute value of $\Psi_{n}$ diverges
to infinity. Therefore, $(\Psi_{n'},\Psi_{n'+1})\in I_{2}\times I_{i}\:(i=1,3)$
is satisfied for a certain $n'$.

\,\,(2-2-b) In the case $(\Psi_{1},\Psi_{2})=(0,0)$, we get $(U_{0},U_{1})=(F,F)$.
Then, we obtain $U_{n}=F$ for all $n$.

(2-3) In the case $(\Psi_{n},\Psi_{n+1})\in I_{2}\times I_{3}$, we
get
\begin{align*}
(U_{n+1},\Psi_{n+1}) & =(F+\Psi_{n},U_{n}-(F+\Psi_{n})),\\
(U_{n+2},\Psi_{n+2}) & =(Q,\Psi_{n}-(Q-F)),
\end{align*}
where
\[
\begin{cases}
\Psi_{n+2}\in I_{1} & (\Psi_{n}\leq(Q-F)-F),\\
\Psi_{n+2}\in I_{2} & (\Psi_{n}>(Q-F)-F).
\end{cases}
\]

\,\,(2-3-a) In the case $\Psi_{n}\leq(Q-F)-F$, we have $(\Psi_{n+1},\Psi_{n+2})\in I_{3}\times I_{1}$.

\,\,(2-3-b) In the case $\Psi_{n}>(Q-F)-F$, we get
\begin{align*}
(U_{n+3},\Psi_{n+3}) & =(F+\Psi_{n}-(Q-F),-\Psi_{n}+2(Q-F)),\\
(U_{n+4},\Psi_{n+4}) & =(Q,\Psi_{n}-2(Q-F)),
\end{align*}
where
\[
\begin{cases}
\Psi_{n+4}\in I_{1} & (\Psi_{n}\leq2(Q-F)-F),\\
\Psi_{n+4}\in I_{2} & (\Psi_{n}>2(Q-F)-F).
\end{cases}
\]

\,\,(2-3-c) In the case $\Psi_{n}\leq2(Q-F)-F$, we have $(\Psi_{n+3},\Psi_{n+4})\in I_{3}\times I_{1}$.

\,\,(2-3-d) In the case $\Psi_{n}>2(Q-F)-F$, repeat the above operation.

If the above discussion did not stop, $\Psi_{n+i}\notin I_{1}$ would
be satisfied all the time. It means $\Psi_{n}<m(Q-F)-F\:\:(m\geq1)$,
which contradicts $-F<\Psi_{n}<Q-F$. Therefore, $(\Psi_{n'},\Psi_{n'+1})\in I_{3}\times I_{1}$
is satisfied for a certain $n'$.

(3-1) In the case $(\Psi_{n},\Psi_{n+1})\in I_{3}\times I_{1}$, we
get
\begin{align*}
(U_{n+1},\Psi_{n+1}) & =(Q,U_{n}-Q),\\
(U_{n+2},\Psi_{n+2}) & =(0,Q),\\
(U_{n+3},\Psi_{n+3}) & =(Q,-Q),\\
(U_{n+4},\Psi_{n+4}) & =(0,Q),\\
 & \vdots
\end{align*}
Therefore, $U_{n}$ alternates between $0$ and $Q$ from a certain
time.

(3-2) In the case $(\Psi_{n},\Psi_{n+1})\in I_{3}\times I_{2}$, it
follows that $(\Psi_{n+1},\Psi_{n+2})\in I_{2}\times I_{i}$ for some
$i$.

(3-3) In the case $(\Psi_{n},\Psi_{n+1})\in I_{3}\times I_{3}$, we
get
\begin{align*}
(U_{n+1},\Psi_{n+1}) & =(Q,U_{n}-Q),\\
(U_{n+2},\Psi_{n+2}) & =(Q,0),
\end{align*}
and thus $(\Psi_{n+1},\Psi_{n+2})\in I_{3}\times I_{2}$.

From the above, a solution of the ultradiscrete equation (\ref{eq:-12})
settles to the periodic solution $(0,Q,0,Q,\ldots)$ except for the
stationary solution $F$. It is can be considered that the periodic
solution $(0,Q,0,Q,\ldots)$ corresponds to the limit cycle (oscillation)
in the original differential equations. The stationary solution $F$
is unstable equilibrium points.

\section*{Acknowledgment}

The author thanks Prof. Seiji Nishioka of Yamagata University for
his invaluable advice.


\begin{thebibliography}{1}
\bibitem{key-1}M. Murata ``Tropical discretization: ultradiscrete
Fisher-KPP equation and ultradiscrete Allen-cahn equation,'' J. differ.
Equ. Appl., 19 (2013), 1008--1021.

\bibitem{key-2}K. Matsuya and M. Murata ``Spatial pattern of discrete
and ultradiscrete Gray-Scott model,'' Discrete Contin. Dyn. Syst.
Ser. B, 20 (2015), 173--187.

\bibitem{key-3}D. Takahashi, A. Shida and M. Usami ``On the pattern
formation mechanism of (2+1)D max-plus models,'' J. Phys. A: Math.
Gen., 34 (2001), 10715--10726.

\bibitem{key-4}R. Hirota and D. Takahashi ``Difference and ultradiscrete
systems,'' Kyoritsu Shuppan, 2015.

\bibitem{key-5}J. H. Merkin ``On wave trains arising in the two-variable
Oregonator model for the BZ reaction,'' IMA. J. Appl. Math., 78 (2013),
513--536.
\end{thebibliography}
\end{document}